\documentclass[a4paper,english]{article}

\usepackage{ifthen}
\usepackage[english]{babel}
\usepackage[utf8]{inputenc}
\usepackage[hyphens]{url}
\usepackage{hyperref}
\usepackage{caption}
\usepackage{amsmath}
\usepackage{amsfonts}
\usepackage{amssymb}
\usepackage{mathtools}
\usepackage[linesnumbered, ruled, noend]{algorithm2e}
\usepackage{graphicx}
\usepackage{enumerate}
\usepackage{xcolor}
\usepackage{xspace}
\usepackage{times}
\usepackage{pgfplots} 
\usepackage[inline]{enumitem}
\usepackage{subfigure}
\usepackage[export]{adjustbox}

\usepackage{tikz}
\usetikzlibrary{patterns}
\usetikzlibrary{positioning}

\usepackage{zibtitlepage}
\usepackage{authblk}

\begin{document}

\title{Warm-starting Strategies in Scalarization Methods for Multi-Objective Optimization}

\author[1]{Stephanie Riedmüller\footnote{corresponding author, riedmueller@zib.de, \ZTPOrcid{0009-0006-4508-4262}}}
\author[1]{Janina Zittel\footnote{zittel@zib.de, \ZTPOrcid{0000-0002-0731-0314}}}
\author[1,2]{Thorsten Koch\footnote{koch@zib.de, \ZTPOrcid{0000-0002-1967-0077}}}

\affil[1]{Zuse Institute Berlin, Berlin, Germany}
\affil[2]{Technische Universtität Berlin, Berlin, Germany}

\maketitle

\begin{abstract}
We explore how warm-starting strategies can be integrated into scalarization-based approaches for multi-objective optimization in (mixed) integer linear programming. Scalarization methods remain widely used classical techniques to compute Pareto-optimal solutions in applied settings. They are favored due to their algorithmic simplicity and broad applicability across continuous and integer programs with an arbitrary number of objectives. 
While warm-starting has been applied in this context before, a systematic methodology and analysis remain lacking. 
We address this gap by providing a theoretical characterization of warm-starting within scalarization methods, focusing on the sequencing of subproblems. However, optimizing the order of subproblems to maximize warm-start efficiency may conflict with alternative criteria, such as early identification of infeasible regions. We quantify these trade-offs through an extensive computational study.
\end{abstract}

\section{Introduction}

Accelerating multi-objective optimization (MOO) has become essential to support timely and effective decision-making for real-world problems.
Given a set $X$ and objective functions $f_1(x), \ldots, f_p(x): X \rightarrow \mathbb{R}$, we consider the MOO problem
\begin{equation}
    \min_{x \in X} (f_1(x), \ldots, f_p(x)). \label{eq:moo}
\end{equation} 
For the definition of \emph{(weakly) efficient} or \emph{Pareto-optimal} solutions in the \emph{decision space} $X$ and the concept of \emph{non-dominance} in the \emph{objective space} $f(X)$, we refer to the standard work \cite{Ehrgott2005}. 

Multi-objective optimization approaches are classified into objective space and decision space algorithms. Scalarization methods, operating in the objective space, remain the most widely used due to their simplicity and broad applicability. Common in practice, they also serve as basis for more advanced algorithms. 
Scalarization converts a multi-objective problem into a set of multiple single-objective problems.

We revisit warm-starting techniques for scalarization methods, i.e., initializing the solving process with information from a previous solver call to improve its performance.
Warm-starting has been applied to scalarization methods in the bi-objective road design problems \cite{akhmet2022} and augmented weighted Tchebycheff network programs \cite{Sun05}. It has also been used in recursive algorithms for multi-objective integer problems \cite{ozlen2014}, in computing enclosures for continuous convex and mixed-integer quadratic problems \cite{eichfelder2023}, and in primal-dual interior point methods for convex quadratic multi-objective problems \cite{molz2006}. 
 
Inspired by these previous applications, we provide a systematic methodology and computational analysis of warm-starting within scalarization methods, focusing on the sequencing of subproblems.

\section{Warm-starting for Scalarization Methods}
We briefly review warm-starting for (mixed-integer) linear programming ((MI)LP), describe the weighted sum and $\varepsilon$-constraint methods, and examine how warm-starting can be applied in each.

\subsection{Warm-starting LPs/MILPs}

A primal or dual feasible solution can be used as as warm-start solution for the process of solving a linear program via the primal or dual simplex algorithm, respectively. In case of a MILP, this supports solving the root relaxation.
A primal feasible solution to a MILP can also be used as a bound in the branch-and-bound tree. 
In general, providing a warm-start solution is recommended, if the first LP-solution takes long or is of bad quality. However, there is no guarantee, that warm-starting will accelerate the solution process.

Note the following observation about primal and dual feasibility depending on the objective function and right-hand side vector of a problem:
Let $x \in \mathbb{R}^n$, $A \in \mathbb{R}^{m \times n}$, $b \in \mathbb{R}^m$ and $c \in \mathbb{R}^n$. Consider the optimization problem
\begin{equation}
    \min c^T x \qquad \text{s.t. } \qquad Ax \geq b, \quad x\geq 0 \label{eq:primal}
\end{equation}
and its dual
\begin{equation}
    \max b^T y \qquad \text{s.t. } \qquad A^Ty \leq c, \quad y\geq 0. \label{eq:dual}
\end{equation}
If $\hat{x}$ is a solution to \eqref{eq:primal} for a given cost vector $c$, then $\hat{x}$ is also primal feasible for any cost vector $c' \neq c$, but is not guaranteed to be dual feasible. 
If $\hat{x}$ is a solution to \eqref{eq:primal} for a given right-hand side vector $b$, then $\hat{x}$ is also dual feasible for any right-hand side vector $b' \neq b$, but is not guaranteed to be primal feasible.

\subsection{Warm-starting the weighted sum method}

The \emph{weighted sum method (WSM)} is to approach \eqref{eq:moo} by solving the single objective problem consisting of a convex combination of the original objectives \cite{Ehrgott2005}:
\begin{equation}
    \min_{x \in X} \sum_{i=1}^{p} w_i f_i(x) \label{eq:ws}
\end{equation}
for weights $w_i > 0$ with $\sum_{i=1}^{p} w_i = 1.$ If a solution is the result of a weighted sum problem \eqref{eq:ws} it is called \emph{supported}, otherwise \emph{non-supported}. Optimal solutions of the weighted sum problem are Pareto-optimal. If the considered problem is convex, the whole set of Pareto-optimal solutions is supported.

Given a solution $\hat{x}$ to \eqref{eq:ws} for some weight vector $w$, then $\hat{x}$ is also primal feasible to \eqref{eq:ws} for another $w' \neq w$, since only the cost vector has been changed.
However, $\hat{x}$ is not guaranteed to be dual feasible, but the primal feasibility allows for warm-starting the initial LP and the branch-and-bound process for a MILP. 
Choosing solutions with similar objective functions as warm-start solution seems to be a promising approach. This could, for example, be achieved by sorting the weight vectors to gain meaningful order of subproblems.

\subsection{Warm-starting the $\varepsilon$-constraint method}
The \emph{$\varepsilon$-constraint method ($\varepsilon$-CM)} is to approach \eqref{eq:moo} by solving the problem for one objective and transforming the other objectives into constraints \cite{Ehrgott2005}:
\begin{equation}
    \min_{x \in X} f_i(x)  \quad \text{s.t.} \quad f_k(x) \leq \varepsilon_k \quad \text{for all } k\neq i \in \{1, \ldots, p\}
\end{equation}
for $\varepsilon \in \mathbb{R}^p$. 
Each efficient (even unsupported) solution can be found by $\varepsilon$-CM, but a solution is only guaranteed to be weakly efficient. There are many variants and hybrid versions of the $\varepsilon$-constraint method. However, we consider the \emph{augmented} $\varepsilon$-constraint method, which guarantees efficient solutions: 
\begin{equation}
    \min_{x \in X} f_i(x) + \rho \sum_{k=1, k \neq i}^{p} f_k(x)  \quad \text{s.t.} \quad f_k(x) \leq \varepsilon_k \quad \text{for all } k\neq i \in \{1, \ldots, p\} \label{eq:aecm}
\end{equation}
for some small, appropriate $\rho>0$.
\begin{figure}
    \centering
    \subfigure[$\varepsilon_2$, $\varepsilon_3$ increasing]{\scalebox{0.8}{\begin{adjustbox}{trim=0pt 30pt 0pt 0pt, clip}
\begin{tikzpicture}

    \begin{scope}
    \fill[gray!20]
      (0.0,2.2)
      .. controls (0.4,3.5) and (2.0,3.5) .. (2.4,2.7)
      .. controls (4.2,1.5) and (2.7,-1.7) .. (1.5,1.2)
      .. controls (0.6,1.4) and (1.0,1.8) .. (0.8,1.9)
      -- cycle;
    \end{scope}

  \begin{axis}[
    width=5cm, height=5cm,
    xmin=0, xmax=4.2, ymin=0, ymax=4.2,
    axis lines=middle,
    xlabel={$f_2$}, ylabel={$f_3$},
    xtick={1,2,3,4},
    ytick={1,2,3,4},
    xticklabels={,,},
    yticklabels={,,},
    grid=major,
    grid style={gray!50},
    axis line style={black},
    tick style={black},
    minor tick num=0,
    clip=false,
    xlabel style={at={(current axis.right of origin)}, anchor=west},
  ylabel style={at={(current axis.above origin)}, anchor=south},
  ]

    \addplot[thick, gray!70, dashed] coordinates {(1,0) (1,4)};
    \node at (axis cs:1,0) [below] {$\varepsilon_2^1$};
    
    \addplot[thick, gray!70, dashed] coordinates {(2,0) (2,4)};
    \node at (axis cs:2,0) [below] {$\varepsilon_2^2$};

    \addplot[thick, gray!70, dashed] coordinates {(3,0) (3,4)};
    \node at (axis cs:3,0) [below] {$\varepsilon_2^3$};
    
    \addplot[thick, gray!70, dashed] coordinates {(4,0) (4,4)};
    \node at (axis cs:4,0) [below] {$\varepsilon_2^4$};

    \addplot[thick, gray!70, dashed] coordinates {(0,1) (4,1)};
    \node at (axis cs:0,1) [left] {$\varepsilon_3^1$};
    
    \addplot[thick, gray!70, dashed] coordinates {(0,2) (4,2)};
    \node at (axis cs:0,2) [left] {$\varepsilon_3^2$};

    \addplot[thick, gray!70, dashed] coordinates {(0,3) (4,3)};
    \node at (axis cs:0,3) [left] {$\varepsilon_3^3$};

    \addplot[thick, gray!70, dashed] coordinates {(0,4) (4,4)};
    \node at (axis cs:0,4) [left] {$\varepsilon_3^4$};

    \coordinate (p1) at (axis cs:0.9,3.8);
    \fill[black] (p1) circle[radius=1.5pt];
    \coordinate (p2) at (axis cs:1.9,3.9);
    \fill[black] (p2) circle[radius=1.5pt];
    \coordinate (p3) at (axis cs:3.0,3.2);
    \fill[black] (p3) circle[radius=1.5pt];
    \coordinate (p4) at (axis cs:0.9,2.8);
    \fill[black] (p4) circle[radius=1.5pt];
    \coordinate (p5) at (axis cs:1.9,2.9);
    \fill[black] (p5) circle[radius=1.5pt];
    \coordinate (p6) at (axis cs:2.8,2.8);
    \fill[black] (p6) circle[radius=1.5pt];
    \coordinate (p7) at (axis cs:3.8,2.1);
    \fill[black] (p7) circle[radius=1.5pt];
    \coordinate (p9) at (axis cs:1.8,1.8);
    \fill[black] (p9) circle[radius=1.5pt];
    \coordinate (p10) at (axis cs:2.8,1.8);
    \fill[black] (p10) circle[radius=1.5pt];
    \coordinate (p11) at (axis cs:3.9,1.4);
    \fill[black] (p11) circle[radius=1.5pt];
    \coordinate (p12) at (axis cs:2.8,0.9);
    \fill[black] (p12) circle[radius=1.5pt];
    \coordinate (p13) at (axis cs:3.8,0.9);
    \fill[black] (p13) circle[radius=1.5pt];

    \draw[red] (axis cs:0,0) -- (axis cs:1,1);
    \draw[red] (axis cs:1,0) -- (axis cs:0,1);
    \draw[red] (axis cs:1,0) -- (axis cs:2,1);
    \draw[red] (axis cs:1,1) -- (axis cs:2,0);
    \draw[red] (axis cs:0,1) -- (axis cs:1,2);
    \draw[red] (axis cs:0,2) -- (axis cs:1,1);

    \node[font=\tiny, gray, fill=white, opacity=0.8] at (axis cs: 0+0.5, 0+0.5) {1};
    \node[font=\tiny, gray, fill=white, opacity=0.8] at (axis cs: 1+0.5, 0+0.5) {2};
    \node[font=\tiny, gray] at (axis cs: 2+0.5, 0+0.5) {3};
    \node[font=\tiny, gray] at (axis cs: 3+0.5, 0+0.5) {4};
    \node[font=\tiny, gray, fill=white, opacity=0.8] at (axis cs: 0+0.5, 1+0.5) {5};
    \node[font=\tiny, gray] at (axis cs: 1+0.5, 1+0.5) {6};
    \node[font=\tiny, gray] at (axis cs: 2+0.5, 1+0.5) {7};
    \node[font=\tiny, gray] at (axis cs: 3+0.5, 1+0.5) {8};
    \node[font=\tiny, gray] at (axis cs: 0+0.5, 2+0.5) {9};
    \node[font=\tiny, gray] at (axis cs: 1+0.5, 2+0.5) {10};
    \node[font=\tiny, gray] at (axis cs: 2+0.5, 2+0.5) {11};
    \node[font=\tiny, gray] at (axis cs: 3+0.5, 2+0.5) {12};
    \node[font=\tiny, gray] at (axis cs: 0+0.5, 3+0.5) {13};
    \node[font=\tiny, gray] at (axis cs: 1+0.5, 3+0.5) {14};
    \node[font=\tiny, gray] at (axis cs: 2+0.5, 3+0.5) {15};
    \node[font=\tiny, gray] at (axis cs: 3+0.5, 3+0.5) {16};

    \draw[->,-stealth, blue, shorten >=2pt, shorten <=2pt] (p1) -- (p2);
    \draw[->,-stealth, blue, shorten >=2pt, shorten <=2pt] (p2) -- (p3);
    \draw[->,-stealth, blue, shorten >=2pt, shorten <=2pt] (p4) -- (p1);
    \draw[->,-stealth, blue, shorten >=2pt, shorten <=2pt] (p4) -- (p5);
    \draw[->,-stealth, blue, shorten >=2pt, shorten <=2pt] (p5) -- (p6);
    \draw[->,-stealth, blue, shorten >=2pt, shorten <=2pt] (p6) -- (p7);
    \draw[->,-stealth, blue, shorten >=2pt, shorten <=2pt] (p9) -- (p10);
    \draw[->,-stealth, blue, shorten >=2pt, shorten <=2pt] (p10) -- (p11);
    \draw[->,-stealth, blue, shorten >=2pt, shorten <=2pt] (p12) -- (p13);
    \draw[->, -stealth, blue, shorten >=2pt, shorten <=2pt]
  (p3) .. controls +(0.0,0.5) and +(0.5,0.5) .. (p3);

  \end{axis}

\end{tikzpicture}

\end{adjustbox}}}
    \hfill
    \subfigure[$\varepsilon_2$, $\varepsilon_3$ decreasing]{\scalebox{0.8}{\begin{adjustbox}{trim=0pt 30pt 0pt 0pt, clip}
\begin{tikzpicture}
    \begin{scope}
    \fill[gray!20]
      (0.0,2.2)
      .. controls (0.4,3.5) and (2.0,3.5) .. (2.4,2.7)
      .. controls (4.2,1.5) and (2.7,-1.7) .. (1.5,1.2)
      .. controls (0.6,1.4) and (1.0,1.8) .. (0.8,1.9)
      -- cycle;
    \end{scope}

  \begin{axis}[
    width=5cm, height=5cm,
    xmin=0, xmax=4.2, ymin=0, ymax=4.2,
    axis lines=middle,
    xlabel={$f_2$}, ylabel={$f_3$},
    xtick={1,2,3,4},
    ytick={1,2,3,4},
    xticklabels={,,},
    yticklabels={,,},
    grid=major,
    grid style={gray!50},
    axis line style={black},
    tick style={black},
    minor tick num=0,
    clip=false,
    xlabel style={at={(current axis.right of origin)}, anchor=west},
  ylabel style={at={(current axis.above origin)}, anchor=south},
  ]

    \addplot[thick, gray!70, dashed] coordinates {(1,0) (1,4)};
    \node at (axis cs:1,0) [below] {$\varepsilon_2^4$};
    
    \addplot[thick, gray!70, dashed] coordinates {(2,0) (2,4)};
    \node at (axis cs:2,0) [below] {$\varepsilon_2^3$};

    \addplot[thick, gray!70, dashed] coordinates {(3,0) (3,4)};
    \node at (axis cs:3,0) [below] {$\varepsilon_2^2$};
    
    \addplot[thick, gray!70, dashed] coordinates {(4,0) (4,4)};
    \node at (axis cs:4,0) [below] {$\varepsilon_2^1$};

    \addplot[thick, gray!70, dashed] coordinates {(0,1) (4,1)};
    \node at (axis cs:0,1) [left] {$\varepsilon_3^4$};
    
    \addplot[thick, gray!70, dashed] coordinates {(0,2) (4,2)};
    \node at (axis cs:0,2) [left] {$\varepsilon_3^3$};

    \addplot[thick, gray!70, dashed] coordinates {(0,3) (4,3)};
    \node at (axis cs:0,3) [left] {$\varepsilon_3^2$};

    \addplot[thick, gray!70, dashed] coordinates {(0,4) (4,4)};
    \node at (axis cs:0,4) [left] {$\varepsilon_3^1$};

    \coordinate (p1) at (axis cs:0.9,3.8);
    \fill[black] (p1) circle[radius=1.5pt];
    \coordinate (p2) at (axis cs:1.9,3.9);
    \fill[black] (p2) circle[radius=1.5pt];
    \coordinate (p3) at (axis cs:3.0,3.2);
    \fill[black] (p3) circle[radius=1.5pt];
    \coordinate (p4) at (axis cs:0.9,2.8);
    \fill[black] (p4) circle[radius=1.5pt];
    \coordinate (p5) at (axis cs:1.9,2.9);
    \fill[black] (p5) circle[radius=1.5pt];
    \coordinate (p6) at (axis cs:2.8,2.8);
    \fill[black] (p6) circle[radius=1.5pt];
    \coordinate (p7) at (axis cs:3.8,2.1);
    \fill[black] (p7) circle[radius=1.5pt];
    \coordinate (p9) at (axis cs:1.8,1.8);
    \fill[black] (p9) circle[radius=1.5pt];
    \coordinate (p10) at (axis cs:2.8,1.8);
    \fill[black] (p10) circle[radius=1.5pt];
    \coordinate (p11) at (axis cs:3.9,1.4);
    \fill[black] (p11) circle[radius=1.5pt];
    \coordinate (p12) at (axis cs:2.8,0.9);
    \fill[black] (p12) circle[radius=1.5pt];
    \coordinate (p13) at (axis cs:3.8,0.9);
    \fill[black] (p13) circle[radius=1.5pt];

    \draw[red] (axis cs:0,0) -- (axis cs:2,1);
    \draw[red] (axis cs:2,0) -- (axis cs:0,1);
    \draw[red] (axis cs:0,1) -- (axis cs:1,2);
    \draw[red] (axis cs:0,2) -- (axis cs:1,1);

    \node[font=\tiny, gray] at (axis cs: 0+0.5, 0+0.5) {16};
    \node[font=\tiny, gray] at (axis cs: 1+0.5, 0+0.5) {15};
    \node[font=\tiny, gray] at (axis cs: 2+0.5, 0+0.5) {14};
    \node[font=\tiny, gray] at (axis cs: 3+0.5, 0+0.5) {13};
    \node[font=\tiny, gray, fill=white, opacity=0.8] at (axis cs: 0+0.5, 1+0.5) {12};
    \node[font=\tiny, gray] at (axis cs: 1+0.5, 1+0.5) {11};
    \node[font=\tiny, gray] at (axis cs: 2+0.5, 1+0.5) {10};
    \node[font=\tiny, gray] at (axis cs: 3+0.5, 1+0.5) {9};
    \node[font=\tiny, gray] at (axis cs: 0+0.5, 2+0.5) {8};
    \node[font=\tiny, gray] at (axis cs: 1+0.5, 2+0.5) {7};
    \node[font=\tiny, gray] at (axis cs: 2+0.5, 2+0.5) {6};
    \node[font=\tiny, gray] at (axis cs: 3+0.5, 2+0.5) {5};
    \node[font=\tiny, gray] at (axis cs: 0+0.5, 3+0.5) {4};
    \node[font=\tiny, gray] at (axis cs: 1+0.5, 3+0.5) {3};
    \node[font=\tiny, gray] at (axis cs: 2+0.5, 3+0.5) {2};
    \node[font=\tiny, gray] at (axis cs: 3+0.5, 3+0.5) {1};
  \end{axis}
\end{tikzpicture}
\end{adjustbox}}}
    \hfill
    \subfigure[$\varepsilon_2$ decr., $\varepsilon_3$ incr.]{\scalebox{0.8}{\begin{adjustbox}{trim=0pt 30pt 0pt 0pt, clip}
\begin{tikzpicture}

    \begin{scope}
    \fill[gray!20]
      (0.0,2.2)
      .. controls (0.4,3.5) and (2.0,3.5) .. (2.4,2.7)
      .. controls (4.2,1.5) and (2.7,-1.7) .. (1.5,1.2)
      .. controls (0.6,1.4) and (1.0,1.8) .. (0.8,1.9)
      -- cycle;
    \end{scope}

  \begin{axis}[
    width=5cm, height=5cm,
    xmin=0, xmax=4.2, ymin=0, ymax=4.2,
    axis lines=middle,
    xlabel={$f_2$}, ylabel={$f_3$},
    xtick={1,2,3,4},
    ytick={1,2,3,4},
    xticklabels={,,},
    yticklabels={,,},
    grid=major,
    grid style={gray!50},
    axis line style={black},
    tick style={black},
    minor tick num=0,
    clip=false,
    xlabel style={at={(current axis.right of origin)}, anchor=west},
  ylabel style={at={(current axis.above origin)}, anchor=south},
  ]

    \addplot[thick, gray!70, dashed] coordinates {(1,0) (1,4)};
    \node at (axis cs:1,0) [below] {$\varepsilon_2^4$};
    
    \addplot[thick, gray!70, dashed] coordinates {(2,0) (2,4)};
    \node at (axis cs:2,0) [below] {$\varepsilon_2^3$};

    \addplot[thick, gray!70, dashed] coordinates {(3,0) (3,4)};
    \node at (axis cs:3,0) [below] {$\varepsilon_2^2$};
    
    \addplot[thick, gray!70, dashed] coordinates {(4,0) (4,4)};
    \node at (axis cs:4,0) [below] {$\varepsilon_2^1$};

    \addplot[thick, gray!70, dashed] coordinates {(0,1) (4,1)};
    \node at (axis cs:0,1) [left] {$\varepsilon_3^1$};
    
    \addplot[thick, gray!70, dashed] coordinates {(0,2) (4,2)};
    \node at (axis cs:0,2) [left] {$\varepsilon_3^2$};

    \addplot[thick, gray!70, dashed] coordinates {(0,3) (4,3)};
    \node at (axis cs:0,3) [left] {$\varepsilon_3^3$};

    \addplot[thick, gray!70, dashed] coordinates {(0,4) (4,4)};
    \node at (axis cs:0,4) [left] {$\varepsilon_3^4$};

    \coordinate (p1) at (axis cs:0.9,3.8);
    \fill[black] (p1) circle[radius=1.5pt];
    \coordinate (p2) at (axis cs:1.9,3.9);
    \fill[black] (p2) circle[radius=1.5pt];
    \coordinate (p3) at (axis cs:3.0,3.2);
    \fill[black] (p3) circle[radius=1.5pt];
    \coordinate (p4) at (axis cs:0.9,2.8);
    \fill[black] (p4) circle[radius=1.5pt];
    \coordinate (p5) at (axis cs:1.9,2.9);
    \fill[black] (p5) circle[radius=1.5pt];
    \coordinate (p6) at (axis cs:2.8,2.8);
    \fill[black] (p6) circle[radius=1.5pt];
    \coordinate (p7) at (axis cs:3.8,2.1);
    \fill[black] (p7) circle[radius=1.5pt];
    \coordinate (p9) at (axis cs:1.8,1.8);
    \fill[black] (p9) circle[radius=1.5pt];
    \coordinate (p10) at (axis cs:2.8,1.8);
    \fill[black] (p10) circle[radius=1.5pt];
    \coordinate (p11) at (axis cs:3.9,1.4);
    \fill[black] (p11) circle[radius=1.5pt];
    \coordinate (p12) at (axis cs:2.8,0.9);
    \fill[black] (p12) circle[radius=1.5pt];
    \coordinate (p13) at (axis cs:3.8,0.9);
    \fill[black] (p13) circle[radius=1.5pt];

    \draw[red] (axis cs:0,0) -- (axis cs:2,1);
    \draw[red] (axis cs:2,0) -- (axis cs:0,1);
    \draw[red] (axis cs:0,1) -- (axis cs:1,2);
    \draw[red] (axis cs:0,2) -- (axis cs:1,1);

    \node[font=\tiny, gray] at (axis cs: 0+0.5, 0+0.5) {4};
    \node[font=\tiny, gray] at (axis cs: 1+0.5, 0+0.5) {3};
    \node[font=\tiny, gray] at (axis cs: 2+0.5, 0+0.5) {2};
    \node[font=\tiny, gray] at (axis cs: 3+0.5, 0+0.5) {1};
    \node[font=\tiny, gray, fill=white, opacity=0.8] at (axis cs: 0+0.5, 1+0.5) {8};
    \node[font=\tiny, gray] at (axis cs: 1+0.5, 1+0.5) {7};
    \node[font=\tiny, gray] at (axis cs: 2+0.5, 1+0.5) {6};
    \node[font=\tiny, gray] at (axis cs: 3+0.5, 1+0.5) {5};
    \node[font=\tiny, gray] at (axis cs: 0+0.5, 2+0.5) {12};
    \node[font=\tiny, gray] at (axis cs: 1+0.5, 2+0.5) {11};
    \node[font=\tiny, gray] at (axis cs: 2+0.5, 2+0.5) {10};
    \node[font=\tiny, gray] at (axis cs: 3+0.5, 2+0.5) {9};
    \node[font=\tiny, gray] at (axis cs: 0+0.5, 3+0.5) {16};
    \node[font=\tiny, gray] at (axis cs: 1+0.5, 3+0.5) {15};
    \node[font=\tiny, gray] at (axis cs: 2+0.5, 3+0.5) {14};
    \node[font=\tiny, gray] at (axis cs: 3+0.5, 3+0.5) {13};

    \draw[->,-stealth, blue, shorten >=2pt, shorten <=2pt] (p4) -- (p1);
    \draw[->,-stealth, blue, shorten >=2pt, shorten <=2pt] (p5) -- (p2);
    \draw[->,-stealth, blue, shorten >=2pt, shorten <=2pt] (p6) -- (p3);
    \draw[->,-stealth, blue, shorten >=2pt, shorten <=2pt] (p7) -- (p3);
    \draw[->,-stealth, blue, shorten >=2pt, shorten <=2pt] (p9) -- (p5);
    \draw[->,-stealth, blue, shorten >=2pt, shorten <=2pt] (p10) -- (p6);
    \draw[->,-stealth, blue, shorten >=2pt, shorten <=2pt] (p11) -- (p7);
    \draw[->,-stealth, blue, shorten >=2pt, shorten <=2pt] (p12) -- (p10);
    \draw[->,-stealth, blue, shorten >=2pt, shorten <=2pt] (p13) -- (p11);

  \end{axis}

\end{tikzpicture}

\end{adjustbox}}}
    \caption{Subset of options for warm-start and infeasibility detection depending on order of subproblems for a 3-objective problem projected on the $f_2/f_3$-plane. Arrows indicate warm-starts and crosses indicate infeasibility detection.}
    \label{fig:ecm}
\end{figure}
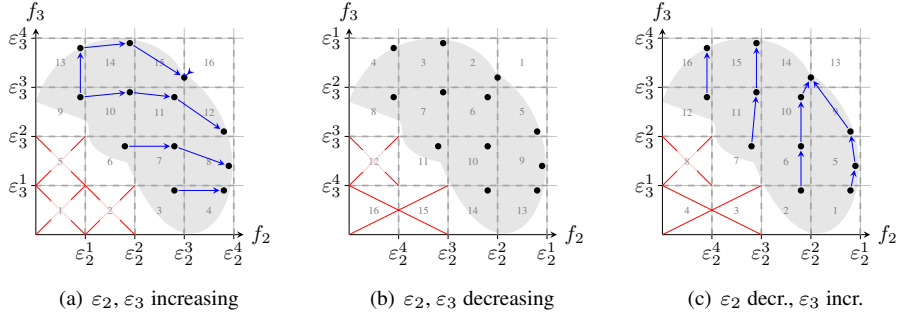

Given a solution $\hat{x}$ to \eqref{eq:aecm} for some constraint vector $\varepsilon$, then $\hat{x}$ is also dual feasible to \eqref{eq:aecm} for another $\varepsilon' \neq \varepsilon$.
It is also primal feasible to \eqref{eq:aecm} for $\varepsilon'$ if $\varepsilon_j \leq \varepsilon'_j$ for $j \in \{1, \ldots, p\}, j \neq i$.
Indeed, it is possible to choose the order of subproblems in a way such that this condition is fulfilled for each subproblem except the first.
While such an order of subproblems is optimal for warm-starting purposes, it stands diametrically against the optimal order for the early detection of infeasible subregions: If problem \eqref{eq:aecm} is infeasible for a given vector $\varepsilon$, then it is also infeasible for $\varepsilon'$ if $\varepsilon_j \geq \varepsilon'_j$ for $j \in \{1, \ldots, p\}, j \neq i$.
Hence, there is a trade-off between the two conflicting targets of warm-starting and detection of infeasibility, see Figure~\ref{fig:ecm}.

\section{Computational study}

We consider 42 test cases of integer problems with 3 and 4 objectives, including assignment problems (AP), knapsack problems (KP), and traveling salesman problems (TSP) \cite{Pettersson2018-3}\cite{Pettersson2018-4}. Additionally, we include two 3-objective mixed-integer unit commitment (UC) problems based on realistic data from district heating networks in a small and a large urban area. All computations are performed on an Intel(R) Xeon(R) CPU E3-1245 v6 @ 3.70GHz. The algorithms are implemented in Python using the commercial solver Gurobi Optimizer version 12.0.1 \cite{gurobi}.
\begin{figure}[h]
    \centering
    \includegraphics[width=0.49\linewidth]{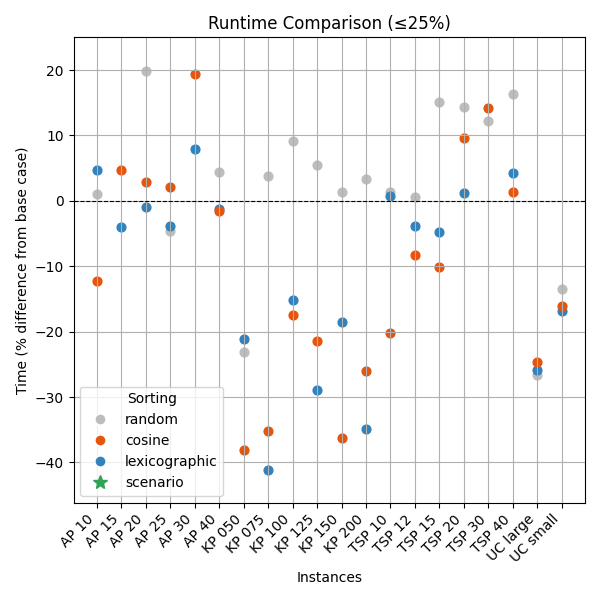}
    \includegraphics[width=0.49\linewidth]{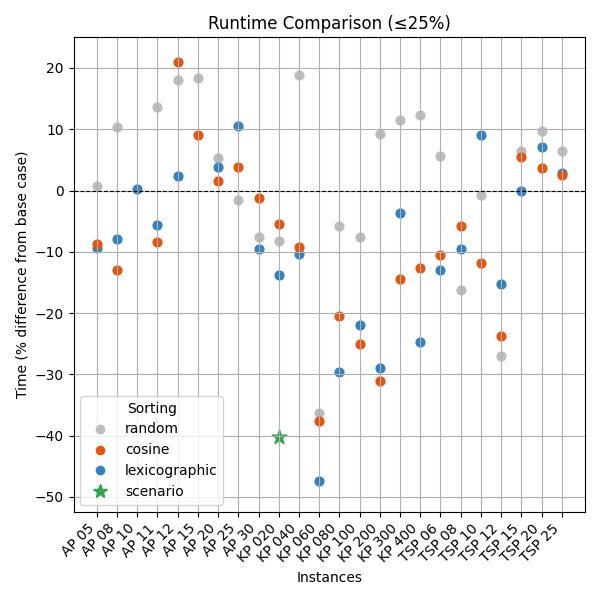}
    \caption{Runtimes for WSM for 3 objectives (left) and 4 objectives (right) in relation to no warmstart. Data points with an increase of runtime of over 25\% have been cut from the diagram.}
    \label{fig:wsm}
\end{figure}

The following setups are investigated:
\begin{itemize}
\item WSM is run with 100 random weight samples. We compare disabled warm-starting to warm-starting from the previous solution under three ordering strategies: random sampling, lexicographic sorting, and sorting by angle against the unit vector.
\item $\varepsilon$-CM is run with 10 equidistant $\varepsilon$-constraints per objective. For each objective, we differentiate, if the order of the $\varepsilon$-constraints is increasing ($+$) or decreasing ($-$) the search region. We further consider three options: \begin{enumerate*}
    \item disabled warm-start and no infeasibility detection
    \item warm-start from the preceding solution and infeasibility detection
    \item warm-start from the previous solution pool and infeasibility detection.
\end{enumerate*}
\item Gurobi supports iterative optimization via Multiple Scenarios\footnote{https://docs.gurobi.com/projects/optimizer/en/current/features/multiscenario.html}, where a base model and its variations are solved in a single run using shared structures. We use this industry-standard to implement WSM and $\varepsilon$-CM: 
the original problem serves as basis, with scenarios varying the objective (WSM) or right-hand sides ($\varepsilon$-CM).
\end{itemize}

Runtimes for warm-started WSM are depicted in Figure~\ref{fig:wsm}. Warm-start effectiveness varies by problem class and is particularly high for KP instances. Sorting weight vectors is typically faster than using random order, reducing runtimes by an average of 27.85\% (3 objectives) and 21.04\% (4 objectives) in the KP instances, though no consistent preference between sorting lexicographically or by angle is observed. Poor-quality warm-starts, such as those from random sorting, even slow down the solver. The Multiple Scenario feature improved the running time only for one instance (KP 020, 4 objectives), otherwise performed over 25\% slower than no warm-start.

\begin{figure}[h]
    \centering
    \includegraphics[width=0.8\linewidth]{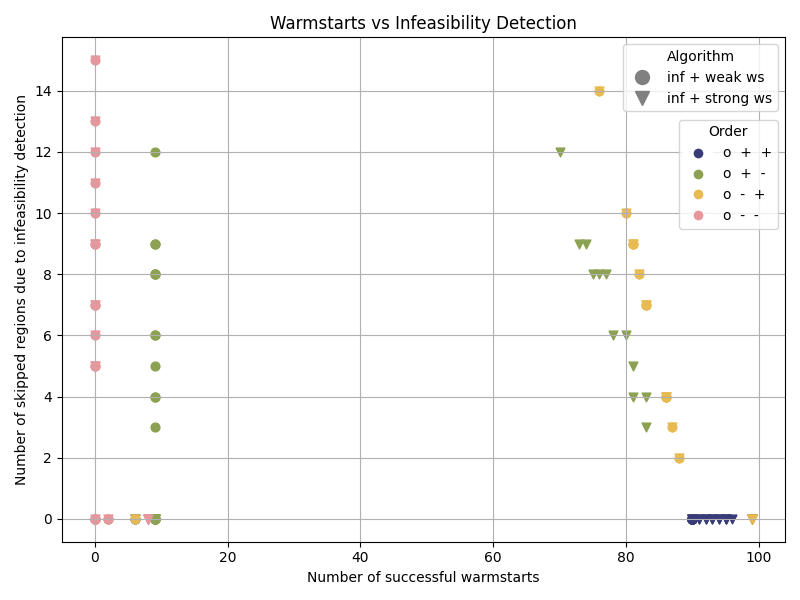}
    \caption{Trade-offs between number of primal feasible warm-starts and detected infeasibilities for 3 objectives in the $\varepsilon$-CM.}
    \label{fig:ec_1}
\end{figure}
For a set of $\varepsilon$-CM variants, Figure~\ref{fig:ec_1} shows the trade-off between the number of primal feasible warm-starts and infeasibility detections. Choosing a subproblem order becomes, therefore, a MOO in itself. Ordering $o - +$ lies entirely on the Pareto front. Ordering $o + +$ and ordering $o - -$ almost exclusively favor warm-starting and infeasibility detection, respectively. Only ordering $o + -$ benefits from using older solutions from the solution pool (strong ws) over just the previous one (weak ws).

\begin{figure}
    \centering
    \includegraphics[width=\linewidth]{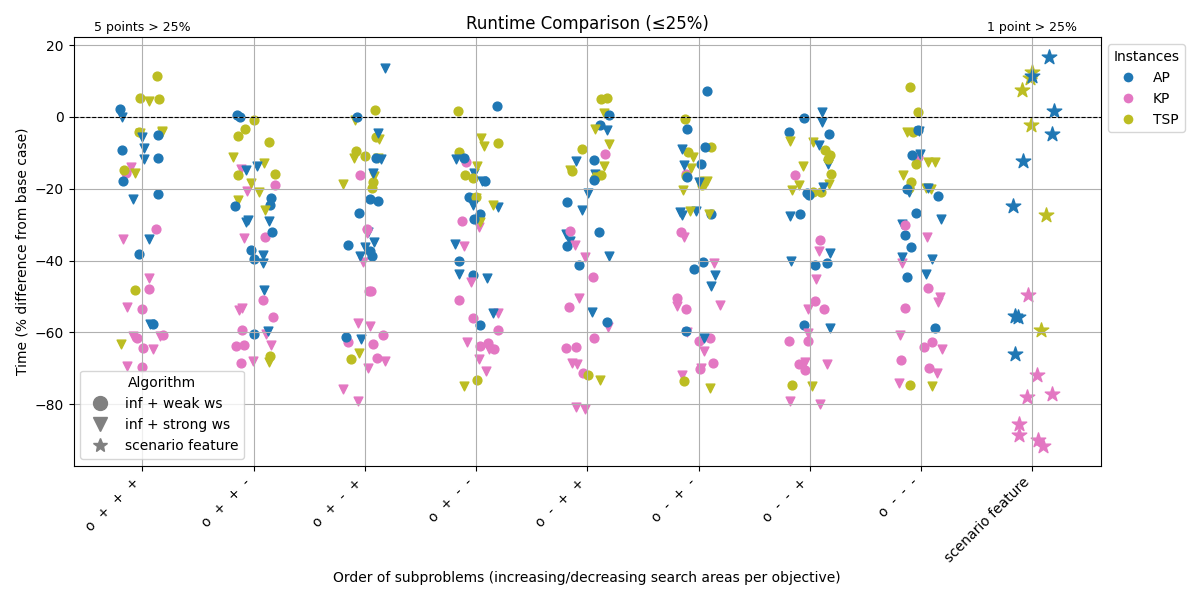}
    \caption{Runtimes for the $\varepsilon$-CM for 4 objectives in relation to no warm-start. Data points with an increase of runtime of over 25\% have been cut from the diagram.}
    \label{fig:ec_2}
\end{figure}

Figure~\ref{fig:ec_2} shows the relative runtime compared to the base version without warm-starts or infeasibility detection, depending on subproblem order. Overall, both techniques reduce runtime by an average of 21.89\% (3 objectives) and 32.39\% (4 objectives). Infeasibility detection tends to impact runtime more than warm-starting when many infeasible regions exist, although the overall effect of the problem order remains small. Strong warm-starts (from earlier solutions) can offer slight improvements over weak ones (from the immediate predecessor), but are costlier to apply (av. improvement: 1.55\% (3 objectives), 3.44\% (4 objectives)). The Multiple Scenario feature performs competitively.

\section{Conclusion}
In general, the impact of warm-starting on scalarization methods depends heavily on the problem structure, making its effectiveness difficult to predict. However, a natural initial solution is available and should be used.
For WSM, we recommend warm-starting by passing the previous solution to the next problem. Sorting weight vectors lexicographically or by angle improves the performance.
For $\varepsilon$-CM, using the last computed solution is also recommended, while searching for better warm-starts offers little added value. Subproblem orders that prevent infeasibility detection should be avoided if many infeasible regions exist. The Multiple Scenario feature has shown to be a good option to implement an out of the box warm-start for the $\varepsilon$-CM.

\section{Acknowledgement}
The work for this article has been conducted in the Research Campus MODAL funded by the Federal Ministry of Research, Technology and Space (BMFTR) (fund numbers 05M14ZAM, 05M20ZBM, 05M2025).

\bibliographystyle{apalike}
\bibliography{IEEEabrv, ref}

\end{document}